\documentclass[12pt]{article}
\usepackage{graphics}
\usepackage{graphicx}
\usepackage{amsmath}
\usepackage{amsfonts}
\usepackage{amssymb}
\everymath{\displaystyle}

\numberwithin{equation}{section}

\begin{document}
\baselineskip 1.5em
\newcommand{\ds}{\displaystyle}

\title{Synthetic Aperture Sonar Imaging
via One-Way Wave Equations}

\author{Quyen Huynh\thanks{Naval Surface Warfare Center - Panama City FL, research partially
supported by the US Office of Naval Research under
N00014-06-WX20559.}
\\
\and Kazufumi Ito\thanks{ Center for Research in Scientific
Computation, Department of Mathematics, North Carolina State
University; research partially supported by the US Office of Naval
Research under N00014-06-1-0067.}}

\maketitle

\begin{abstract}
We develop an efficient algorithm for Synthetic Aperture Sonar
imaging based on the one-way wave equations. The algorithm
utilizes the operator-splitting method to integrate the one-way
wave equations. The well-posedness of the one-way wave equations
and the proposed algorithm is shown. A computational result
against real field data is reported and the resulting image is
enhanced by the BV-like regularization.
\end{abstract}

\section{Introduction}

In this paper we discuss the migration method based on one-way wave
equations \cite{C,G} for Synthetic Aperture Sonar (SAS) imaging
\cite{H}. The one-way wave equation integrates the data within a
given angle and minimizes the undesirable effects of unwanted
reflections. Efficient and stable integration methods of the one-way
wave equation based on the operator splitting method are used to
develop a fully discretized algorithm. The stability analysis and
the required operation count of the proposed algorithm are given. We
test the proposed method for real field data and report our SAS
imaging results.  We also discuss the image enhancement method for
the resulting images, based on BV-like regularization technique
\cite{IK}.

In side-scan (side-looking) sonar systems a platform containing a
moderately large real aperture antenna travels along a rectilinear
path in the along track direction and periodically transmits a pulse
at an angle that is perpendicular to the platform path.
These systems produce strip-map (SAS) images . A strip-map image
is built up as follows; the imaging system operates such that the echoes
from the current pulse are received before the next pulse is transmitted.
As these echoes are received they are demodulated,
pulse compressed, and detected (only the magnitude information is retained).
Each detected pulse produces a range line of the real aperture image.
As the platform moves these range lines are displayed
next to each other at pixel spacings that scale relative to the along track spacing
of the pulses $\Delta x=v_p\tau$ where $v_p$ is the platform velocity and $\tau$
is the pulse repetition period.  The final image
is essentially a raster scan of a strip of the sea floor, hence
the name "strip-map image". Synthetic aperture imaging is a coherent imaging
technique that exploits the extra information
available in the phase of the real aperture data.
We adopt the Stop and Go model;
a point source  radiates at time $t=0$, a spherical wave that
reaches the sampling points after different time intervals.
If the source is placed at $(x_0,z_0)$ the time $t(x,x_0,z_0)$ at
which the wave arrives at the sampling point $(x,z)$ is:
$$
t(x,z,x_0,z_0)=\frac{2}{c} \sqrt{(x-x_0)^2+(z-z_0)^2}.
$$
The field $d$ due to a distribution $s(x,z)$ of source emitting at $t=0$
can be expressed by
$$
\hat{d}(x,z,\omega)=\frac{1}{4\pi}\int s(x^\prime,z^\prime)
\frac{e^{-j(2\omega/c)\sqrt{(x-x^\prime)^2+(z-z^\prime)^2}}}
{\sqrt{(x-x^\prime)^2+(z-z^\prime)^2}}\,dx^\prime dz^\prime
$$
where $\hat{d}$ is the Fourier transform (in time) of the signal
$d$. SAS measures
$$
\mbox{SAS}(x,t)=d(x,z=0,t)
$$
along the sonar path $(x,z=0)=\Gamma$.

Thus, SAS imaging is formulated as a linear inverse problem;
\vspace{2mm}

\noindent\underline{Problem:} Reconstruct $s(x,z)$ from SAS data
SAS$(x,t)$. \vspace{2mm}

Among a number of algorithms \cite{H,HP,So} and reference therein,
which have been developed for Problem the frequency domain
$\omega$-$k$ method based on Stolt's map \cite{H,S,C} is the most
efficient and accurate method. As will be discussed in Section 4 it
has certain limitations, especially it assumes the homogeneous
scattered media. The proposed method can incorporate inhomogeneous
media and has additional capabilities, (see Section 4).

An outline of the paper is as follows. In Section 2 we describe a
geometric migration method based on one-way wave equations for
reconstructing $s$. A noble algorithm using the integration of the
one-way wave equations based on the operator-splitting method is
developed and its stability and complexity are analyzed in Section
3. In Section 3 we list advantages of the proposed method comparing
to the $\omega$-$k$ method. The image enhancement technique based on
the BV-type reguralization is discussed in Section 4. In Section 5
we present a test against real field data, provided by the Naval
Surface Warfare Center-Panama City, Florida and a comparison with
the $\omega$-$k$ method.

\section{Geometric Migration}

We construct an approximating
solution based on the geometrical migration via the one-way wave equations.
Let
$$
A(k_x,\omega)={\cal F}_{x,t}\, \mbox{SAS}(x,t).
$$
Assume the plane wave extrapolation
$$
D(k_x,k_z,\omega)=A(k_x,\omega)exp(j(\omega t+k_x x+k_z z))
$$
with
$$
\omega^2=\frac{c^2}{4}(k_x^2+k_y^2).
$$
Then the inverse Fourier transform of $D$
$$
\tilde{d}(x,z,t)=\frac{1}{(2\pi)^3}\int D(k_x,k_z,\omega) dk_x dk_z d\omega
$$
satisfies the wave equation
$$
\frac{4}{c^2}\frac{\partial^2 \tilde{d}}{\partial t^2}
=\frac{\partial^2 \tilde{d}}{\partial x^2}
+\frac{\partial^2 \tilde{d}}{\partial z^2}
\leqno (3) $$
with the boundary condition at $z=0$
$$
\tilde{d}(x,0,t)=\mbox{SAS}(x,t)
$$
and
$$
\tilde{d}(x,z,T)=0  \mbox{  and   }
\frac{\partial \tilde{d}}{\partial t}(x,z,T)=0.
$$
Wave equation based migration integrates the wave equation (3)
backward in time to obtain an approximation $\tilde{s}$ of
distribution $s$ as;
$$
\tilde{d}(x,z,0)=\tilde{s}(x,d).
$$

SAS data is created by integrating over the beam-width of the
sensor. The radiation pattern of any dimension (width or length)
of an aperture has an angular dependence that is referred to as
the beam pattern of the aperture. Beam patterns are frequency
dependent and have beam-widths given by the 3dB response of their
main lobes; $\ds \theta=\alpha_w\,\frac{c}{fD}$ where $D$ is the
length of the aperture and $f$ are the frequency of
the signal that the aperture is transmitting or receiving.  The
term $\alpha_w$ is a constant reflecting the main lobe widening
due to weighting of the aperture illumination function. For
example $f=120$kHz and $D=0.04$m and $\alpha_w=1$ gives $\theta = 17.9$ degrees
Thus, in order to speed-up the wave equation based algorithm and
minimize the undesirable effects of unwanted reflections we use
the (15 degree) one-way wave equation based on
$$
k_z=k\sqrt{1-\left(\frac{k_x}{k}\right)^2} \sim
k\,(1-\frac{1}{2}(\frac{k_x}{k})^2)
\leqno (4) $$
where $k=2\omega/c$ and we assumed $|k_x/k|<<1$.
In time domain (4) is equivalently written as
$$
\frac{4}{c^2}\frac{\partial^2 u}{\partial t^2}
+\frac{2}{c}\frac{\partial^2 u}{\partial z\partial t}
=\frac{1}{2}\frac{\partial^2 u}{\partial x^2}.
\leqno (5) $$
with
$$
u(t,x,0)=\mbox{SAS}(x,t),\;x\in\Gamma.
$$

An advantage of the method is that it allows one to have a specified variable wave
speed $c=c(x,z)$ of media. The corresponding method for the polar and
cylindrical geometry is given as

\noindent\underline{Polar coordinate}
$$
\frac{4}{c^2}\,\frac{\partial^2 u}{\partial t^2}
+\frac{2}{c}\,\frac{\partial^2 u}{\partial \nu\partial t}
=\frac{1}{2}\frac{1}{r}\,\frac{\partial^2 u}{\partial\theta^2}.
$$

\noindent\underline{Cylinder}
$$
\frac{4}{c^2}\,\frac{\partial^2 u}{\partial t^2}
+\frac{2}{c}\,\frac{\partial^2 u}{\partial \nu\partial t}
=\frac{1}{2}\,(\frac{1}{r}\frac{\partial^2 u}{\partial\theta^2}
+\frac{\partial^2 u}{\partial z^2}).
$$

We can derive the wide angle one-way wave equation based on the
rational approximation
$$
k_z=k \sqrt{1-\left(\frac{k_x}{k}\right)^2} \sim k\,(1-\frac{\alpha (k_x/k)^2}
{1-\beta (k_x/k)^2})
\leqno (6) $$
we have
$$
k_z(k-\frac{\beta}{k} k_x^2)=k^2-(\alpha+\beta) k_x^2
$$
The differential form is given by
$$
\frac{4}{c^2}\, \frac{\partial^2 u}{\partial t^2}
+\frac{\partial}{\partial z}\left(\frac{2}{c}\frac{\partial
u}{\partial t} -\beta\,\frac{c}{2}\int \frac{\partial^2
u}{\partial x^2}\,dt \right)=(\alpha+\beta)\frac{\partial^2
u}{\partial x^2} \leqno (7)
$$
With $\alpha=.5,\;\beta=.25$ and $\alpha=.478,\; \beta=.376$,
(7) is called 45 degree and 65 degree approximation, respectively.

\section{Migration by the operator splitting}

With normalization of the time (t) by the wave speed $\ds
\frac{c}{2}$ and reverting the time, (5) is written as
$$
\left(\begin{array}{c} u_t \\ \\ v_t \end{array}\right)
=\left(\begin{array}{cc} 0 & 0 \\ \\ 0 & -\frac{\partial}{\partial
z}\end{array}
\right) \left(\begin{array}{c} u\\ \\ v \end{array}\right)
+ \left(\begin{array}{cc} 0 & 1 \\ \\ \frac{1}{2}
\frac{\partial^2}{\partial x^2} & 0 \end{array}\right)
\left(\begin{array}{c} u \\ \\ v\end{array} \right).
\leqno (8) $$
So, we apply the time splitting on $[t,t+\Delta t]$ of the Lie-Trotter form
$$
\left(\begin{array}{c} u_t\\ \\ v_t \end{array}\right)
=\left(\begin{array}{cc} 0 & 0 \\ \\ 0 & -\frac{\partial}{\partial
z}\end{array} \right)\left(\begin{array}{c} u\\ \\ v
\end{array}\right),\quad \left(\begin{array}{c} u_t\\ \\ v_t
\end{array} \right) = \left(\begin{array}{cc} 0 & 1 \\ \\
\frac{1}{2} \frac{\partial^2}{\partial x^2} & 0 \end{array}\right)
\left(\begin{array}{c} u \\ \\ v \end{array} \right). \leqno (9)
$$
The first step of (9) is equivalent to the shift operation;
$$
\left\{ \begin{array}{l} v(t+\Delta t,x,z)=v(t,x,z-\Delta t)
,\quad z \ge \Delta t
\\ \\
v(t+\Delta t,x,z)=\frac{\partial}{\partial t}\mbox{SAS}(x,t+z),
\quad z\in [0,\Delta t) \end{array} \right.
$$
The second step of (9) is the one-D wave equation in $x$ and is
well-posed. In fact, let $\Omega=[-L,L]\times [0,1]$ and
$H^{1,x}(\Omega)=\{\phi\in L^2(\Omega): \frac{\partial}{\partial
x}\phi \in L^2(\Omega)\}$. Let $X_1=H^{1,x}(\Omega) \times
L^2(\Omega)$ be the Hilbert space equipped with
$$
|(u,v)|_{X_1}^2=\int_\Omega (|\frac{\partial u}{\partial
x}|^2+2|v|)^2)\, dxdz.
$$
Define the linear operator ${\cal A}_1$ on $X_1$ by
$$
{\cal A}_1(u,v)=(v,\frac{1}{2}\frac{\partial^2 u}{\partial x^2})
$$
with
$$\begin{array}{l}
\ds dom({\cal A}_1)=\{(u,v)\in X_1: v \in H^{1,x}(\Omega),\;
\frac{\partial^2 u}{\partial x^2}\in L^2(\Omega)
\\ \\
\ds\qquad\qquad  \mbox{  with  }
\frac{\partial u}{\partial x}(\pm L,z)=0\}
\end{array} $$
Then, ${\cal A}_1$ is dissipative and skew-adjoint on $X_1$ and
thus generates a strongly continuous group on $X_1$. Hence, it is
easy to show that if $(u,v)$ is generated by (9) then
$$
|(u,v)(t+\Delta t)|_{X_1}^2 \le  |(u,v)(t)|^2_{X_1}
+\int^{t+\Delta t}_t |\frac{\partial}{\partial
t}\mbox{SAS}(x,s)|^2\,dxds
$$
and
$$
|(u,v)(T)|_{X_1}^2 \le \int^T_0 |\frac{\partial}{\partial
t}\mbox{SAS}(x,t)|^2\,dxdt.
$$
Similarly, we can argue that (8) itself is well-posed, i.e., if we
define the operator ${\cal A}$ on $X_1$ by
$$
{\cal A}(u,v)=(v,div_{x,z}(\frac{1}{2}\frac{\partial u}{\partial
x},-v))
$$
with
$$ dom({\cal A})=\{(u,v)\in X_1:v\in H^{1,x}(\Omega),
\;div_{x,z}(\frac{1}{2}\frac{\partial u}{\partial x},-v) \in
L^2(\Omega) \mbox{  with  } v(x,z)=0,\;\frac{\partial u}{\partial
x}(\pm L,z)=0\},
$$
then ${\cal A}$ is dissipative and generates a contractive,
strongly continuous semigroup on $X_1$.

We fully discretize (9) and obtain \vspace{2mm}

\noindent\underline{\bf Algorithm I}
$$\begin{array}{l}
\hat{v}^{n+1}_{i,j+1}=v^n_{i,j},\;\; 1 \le j \le n \mbox{  with  }
\hat{v}^{n+1}_{i,0}=\ds
\frac{\mbox{SAS}^{n+1}_i-\mbox{SAS}^n_i}{\Delta t}
\\ \\
u^{n+1}_{\cdot,j}=(I+\frac{\tilde{c}^2}{2}H)^{-1}
(u^n_{\cdot,j}+\Delta t \hat{v}^{n+1}_{\cdot,j}),\;\;\;
v^{n+1}_{\cdot,j}=\frac{u^{n+1}_{\cdot,j}-u^n_{\cdot,j}}{\Delta t},
\; 1\le j \le \min(n,M)
\end{array} \leqno(10) $$
where $u^n_{i,j}$ and $v^n_{i,j}$ represents the value
of $u$ and $v$ at the grid-point $(i\Delta x,j\Delta z)$ at time $n\Delta t$,
respectively. Here, $\Delta t=\Delta z$ and  $\ds \tilde{c}
=\frac{\Delta z}{\Delta x}$, $H \in R^{N+1,N+1}$ is the tri-diagonal matrix
defined by
$$\begin{array}{l}
(Hu)_i=-(u_{i+1}-2u_i+u_{i-1}),\;\; 2 \le i \le N,
\\ \\
\quad\mbox{and}\quad
(Hu)_1=-(u_2-u_1),\quad (Hu)_{N+1}=u_{N+1}-u_N
\end{array} $$
and corresponds to the central difference approximation of $\ds
-\frac{\partial^2 u}{\partial x^2}$. Also, we used the implicit
Euler scheme to  integrate the second step (1-D wave equation in $x$).
That is,
$$
\frac{u^{n+1}_i-u^n_i}{\Delta t}=-\frac{1}{\Delta x^2}(Hu)_i,
\quad \frac{v^{n+1}-\tilde{v}^n}{\Delta t}=u^{n+1}.
$$
The number of operations at the n-th time step of (10) is of order
$O(N\,\min(n,M))$. $M$ is the number of the focusing step at each
pixel $(i,j)$ in cross-range direction $x$ and if $j \ge M$, then
$u^{n+1}_{i,j+1}=u^n_{i,j}$. Thus, the total operation is of order
$O(MM^2)$.

For the wide angle equation (7) we define
$$
F=\frac{2}{c}\frac{\partial u}{\partial t}-\beta\frac{c}{2} \int^t_0
\frac{\partial^2 u}{\partial x^2}\,dt \quad \mbox{and} \quad
v=\frac{2}{c}\frac{\partial u}{\partial t}.
$$
It follows from (7) that
$$
\frac{2}{c}\frac{\partial F}{\partial t}
=(\frac{2}{c})^2\frac{\partial^2 u}{\partial t^2}
-\beta\frac{\partial^2 u}{\partial x^2}=
-\frac{\partial F}{\partial z}+\alpha\,\frac{\partial^2 u}{\partial x^2}
$$
and
$$
\frac{\partial}{\partial t}(v-F)=
\frac{c}{2}\beta\,\frac{\partial^2 u}{\partial x^2}.
$$
Thus, (7) is equivalent to
$$\begin{array}{l}
\ds \frac{2}{c}\frac{\partial F}{\partial t}+\frac{\partial
F}{\partial z}= \alpha\,\frac{\partial^2 u}{\partial x^2}
\\ \\
\ds \frac{2}{c} \frac{\partial u}{\partial t}=(v-F)+F
\\ \\
\ds \frac{2}{c}\frac{\partial}{\partial t}(v-F)=\beta\,
\frac{\partial^2 u}{\partial x^2}.
\end{array} \leqno (11) $$
With $\tilde{v}=v-F$, we use the three step splitting:
$$
\left\{ \begin{array}{l}
\ds \frac{\partial F}{\partial t}+\frac{\partial F}{\partial z}=0
\\ \\
\ds \frac{\partial u}{\partial t}=0
\\ \\
\ds \frac{\partial \tilde{v}}{\partial t}=0
\end{array} \right.\quad
\left\{ \begin{array}{l}
\ds \frac{\partial F}{\partial t}=0
\\ \\
\ds \frac{\partial u}{\partial t}=\tilde{v}
\\ \\
\ds \frac{\partial \tilde{v}}{\partial t}=\beta\,\frac{\partial^2 \tilde{v}}
{\partial x^2} \end{array} \right. \quad
\left\{ \begin{array}{l}
\ds \frac{\partial F}{\partial t}=\alpha\,\frac{\partial^2 F}
{\partial x^2}
\\ \\
\ds \frac{\partial u}{\partial t}=F
\\ \\
\frac{\partial \tilde{v}}{\partial t}=0
\end{array} \right.
\leqno (12) $$ If $\beta=0$ then $\tilde{v}=0$, $F=v$ and thus it
reduces to the two-step splitting method (9). The first equation
is accompanied by the boundary condition
$$
F(t,x,0)=\frac{\partial}{\partial t}\mbox{SAS}(t,x)-\tilde{v}(t,x,0).
$$
Each step of (12) is a well-posed linear system as shown above and
we can prove that (11) is well-posed. In fact, let
$\Omega=[-L,L]\times [0,1]$ and define the linear operator on
${\cal A}_2$ on $X_2=L^2(\Omega) \times H^{1,x}(\Omega) \times
L^2(\Omega)$ by
$$
{\cal A}_2(F,u,\tilde{v})=(-\frac{\partial F}{\partial z}+\alpha
\frac{\partial^2 u}{\partial x^2}, \tilde{v}+F, \beta
\frac{\partial^2 u}{\partial x^2})
$$
with
$$\begin{array}{l}
\ds dom({\cal A}_2)=\{\frac{\partial}{\partial z}F \in
L^2(\Omega), \mbox{  with  } F(\cdot,0)=0 \mbox{  and}
\\ \\
\ds\qquad
\frac{\partial^2}{\partial x^2}u \in L^2(\Omega)
\mbox{  with  }  \frac{\partial u}{\partial x}(\pm L,z)=0,
\frac{\partial}{\partial x}(\tilde{v}+F)\in L^2(\Omega)\}.
\end{array} $$
We equip $X_2$ with norm
$$
|(F,u,\tilde{v})|^2_{X_2}=\int_\Omega (|\frac{\partial}{\partial
x}u|^2 +\frac{1}{\alpha}|F|^2+\frac{1}{\beta}|\tilde{v}|^2)\,dxdz
$$
Then, ${\cal A}_2$ is dissipative, i.e.,
$$\begin{array}{l}
({\cal A}_2(F,u,\tilde{v}),(F,u,\tilde{v}))
\\ \\
\ds\quad =
\int_\Omega (\frac{\partial^2 u}{\partial x^2}(\tilde{v}+F)
+\frac{\partial u}{\partial x} \frac{\partial}{\partial x}(\tilde{v}+F)
-\frac{\partial F}{\partial z}F)\,dxdz
\\ \\
\ds \quad
=-\frac{1}{2}\int^{L}_{-L} |F(x,1)|^2\,dx \le 0.
\end{array} $$
Since $range({\cal A}_2)=X_2$, ${\cal A}_2$ generates a  strongly
continuous, contraction semigroup on $X_2$. Similarly, we have the
energy estimate
$$
\int (|\frac{\partial}{\partial x}u(T)|^2
+\frac{1}{\alpha}|F(T)|^2+\frac{1}{\beta}|\tilde{v(T)}|^2\,dxdz
\le \int^T_0\int|F(t,x,0)|^2\,dx\,dt.
$$
Algorithm I is extended to integrate (12) as follows; \vspace{2mm}

\noindent\underline{\bf Algorithm II}
$$\begin{array}{l}
\hat{F}^{n+1}_{i,j+1}=F^n_{i,j},\;\; 1 \le j \le n \mbox{  with  }
\hat{F}^{n+1}_{i,0}=\ds
\frac{\mbox{SAS}^{n+1}_i-\mbox{SAS}^n_i}{\Delta t}-
\tilde{v}^n_{i,0}
\\ \\
\hat{u}^{n+1}_{\cdot,j}=(I+\beta\tilde{c}^2 H)^{-1}
(u^n_{\cdot,j}+\Delta t \hat{F}^{n+1}_{\cdot,j}),\;\;\;
F^{n+1}_{\cdot,j}=\frac{\hat{u}^{n+1}_{\cdot,j}-u^n_{\cdot,j}}{\Delta t},
\; 1\le j \le \min(n,M)
\\ \\
u^{n+1}_{\cdot,j}=(I+\alpha\tilde{c}^2 H)^{-1}
(\hat{u}^{n+1}_{\cdot,j}+\Delta t \tilde{v}^n_{\cdot,j}),\;\;\;
v^{n+1}_{\cdot,j}=\frac{u^{n+1}_{\cdot,j}-\hat{u}^{n+1}_{\cdot,j}}{\Delta t},
\; 1\le j \le \min(n,M).
\end{array} $$
That is, we require double the operations for the integration of
Algorithm II.

\section{Advantages of the proposed methods}

The frequency domain $\omega$-$k$ method based on Stolt's map
\cite{S} is the most efficient and accurate method for the
homogeneous media due to the efficiency of fast Fourier transform.
It also assumes a rectilinear sonar path.

We can use our proposed algorithms as a means to compensate the
motion of sonar path. That is, let $\Gamma$ be a curved sonar path
and $\Gamma_0$ is a reference rectilinear path (z=0). Then we solve
(5) or (7) on the domain enclosed by the boundaries $\Gamma$ and
$\Gamma_0$ with boundary value
$$
u(t,x,z)=\mbox{SAS}(t,x),\quad (x,z) \in \Gamma
$$
In this way we have the mapped-SAS data $u(t,x,0)$ at $\Gamma_0$
and then apply the $omega-k$ method for the rectangular domain $\Omega$.

Our implementation (10) of the one-way wave equations is easily adjusted
to the case of layered media $c=c(z)$ by varying the range increments
$\Delta z$.

The proposed method can allow to localize the integration on
sub-layered regions (assuming the homogeneous media). Also, we can
integrate (5) or (7) in overlapped sub-domains in the cross-range
(x) direction and then apply the superposition. This improves the
efficiency of the proposed algorithms.

\section{BV-type Regularization for Enhancement of SAS imaging}

SAS imaging $s(x,z)$ may be altered by inhomogeneity of the field,
sensor noise and irregularity of the sonar path and so on. We use
the image enhancement technique based on BV-type reguralization
\cite{IK}.

\noindent\underline{Enhancement $S$} of $s$ minimizes
$$
\int_\Omega |S-s|^2\,dxdz +\beta \int_\Omega \varphi(
|\frac{\partial S}{\partial x}|^2 +|\frac{\partial S}{\partial
z}|^2)\,dx dz \leqno (12) $$
where
$$
\mbox{$\beta>0$ is the regularization parameter}
$$
and
$$
Q(\phi)=\int_\Omega \varphi(|\nabla S|^2)\,dx dz \mbox{  defines
the restoration energy}.
$$
The followings summarize our findings in \cite{IK} on the
enhancement based on(12);

\begin{itemize}

\item $\varphi(t^2)=t^2$ corresponds to the standard Gaussian filter
and works well for a smooth image $s$.

\item $\varphi(t^2)=t$ corresponds to the BV (nonlinear) filter
and restores edges and flat regions of image $s$ very well.
But, it has significant stair-case effects.

\item In order to deal with images with multi-scales
of edges, flat, and smooth regions
we developed an algorithm which uses
$$
\varphi^\prime(s)=\left\{
\begin{array}{ll} \frac{1}{\sqrt{s}} & s\in [1,\infty) \\
1 & s\in [\delta,1]
\\ \frac{1}{\sqrt{s}} & s\in (0,\delta)
\end{array} \right.
$$
It is based on the scale analysis and we demonstrated the
applicability of the algorithm in \cite{IK}.

\item The necessary and sufficient condition of (12) is given by
$$
-\beta\,\nabla\cdot(\varphi^\prime(|\nabla S|^2)\nabla S)+S=s.
$$
\end{itemize}

An efficient algorithm for finding $S$ based on the fixed point iterate;
$$
-\beta\,\nabla\cdot(\varphi^\prime(|\nabla S^k|^2)\nabla S^{k+1})+S^{k+1}=s
$$
is developed and analyzed in \cite{IK} and is used in our test.

\section{A Test}

The algorithm is successfully applied to real data that are
available to us via the Naval Surface Warfare Center (NSWC) and
shows a promising capability.  A full capability is going to be
tested in the line of its advantages discussed in Section 4. In a
CRSC tereport, CRSC-TR09-12 at
http://www.ncsu.edu/crsc/reports/reports09.htm we show the raw SAS
data, SAS imaging by algorithm (10), and the image enhanced by our
enhancement algorithm.






\end{document}